\documentclass[11pt]{article}
\usepackage[T1]{fontenc}
\usepackage{lmodern}
\usepackage{amsmath,amssymb,amsthm,mathtools}
\usepackage{booktabs}
\usepackage{float}
\usepackage{array}
\usepackage{enumitem}
\usepackage{xcolor}
\usepackage{xurl}
\usepackage[colorlinks=true,linkcolor=blue,citecolor=blue,urlcolor=blue]{hyperref}
\usepackage[a4paper,margin=1in]{geometry}
\newcommand{\mathlib}{\texttt{Mathlib}}
\newcommand{\code}[1]{\texttt{#1}}
\newcommand{\codeurl}[1]{\nolinkurl{#1}}

\title{A Lean 4 Verification Report for\\
\emph{Subregular Affine Cells and the Level $-1$ Vertex Algebra of Type $D$}}
\author{Sihai Jin\\
Department of Mathematics, Sichuan University\\
\texttt{jinsihai@stu.scu.edu.cn}}
\date{}

\begin{document}
\maketitle

\begin{abstract}
We report a Lean 4 formal verification accompanying the paper
\emph{Subregular affine cells and the level $-1$ vertex algebra of type $D$}
(arXiv:2608.11997).  The formalization kernel-checks substantial internal
parts of the proof architecture, including the Section~4 membership/descent
chain, the type-$D$ norm-gap argument, zero-orbit energy and
signed-permutation rigidity, node-weight and numerical rigidity calculations,
the exhaustion logic, the passage to the candidate quotient classification,
the simple-object count, the final additive Grothendieck-group
comparison, and the coefficient-substitution layer of the uniform character
formula.  Higher representation-theoretic results
whose foundational infrastructure is not presently constructed in the file
are isolated as explicit semantic interfaces rather than introduced as Lean
axioms.  Thus the precise claim is a kernel-checked internal deduction from
explicit representation-theoretic boundary inputs, not a from-scratch
formalization of vertex algebras, BRST reduction, finite $W$-algebras, or
affine Hecke theory inside \mathlib.
\end{abstract}

\section{Purpose and relation to the original paper}
The original paper proves the Shan--Yan--Zhao affine left-cell conjecture for
$L_{-1}(D_\ell)$, $\ell\ge 5$.  Its main theorem classifies the simple objects
in the distinguished vacuum block and identifies the Grothendieck group with
the $q=1$ specialization of the corresponding dual affine left-cell module
\cite{Jin2026}.

The accompanying Lean source is a single-file verification bundle.  Its
purpose is not to recreate all background representation theory from basic
definitions.  Instead, it separates the proof into two layers:
\begin{enumerate}[label=(\roman*)]
  \item calculations and logical deductions checked directly by the Lean
  kernel; and
  \item external representation-theoretic facts represented by explicit
  hypotheses or structured interfaces.
\end{enumerate}
The final top-level theorem is
\codeurl{DTypeMainTheorem.theorem_1_1}.  The released source is available both
as arXiv ancillary material and in the public GitHub repository
\url{https://github.com/jshemail12345-debug/SubregularAffineCells-Lean}; the
version corresponding to this report is frozen as release
\url{https://github.com/jshemail12345-debug/SubregularAffineCells-Lean/releases/tag/v1.0.0}.

\section{Scope of the formalization}
\subsection{Kernel-checked internal components}
The single Lean file contains the following principal internally verified
components.
\begin{enumerate}
  \item The concrete rank-one parameter identities in Lemma~4.8, including
  $c_\beta(-m-2)=c_\beta(m)$, the special parameter $m=-2$, and the
  scalar specialization of Premet's second parameter at $m=0$.
  \item The Appendix-A coordinate comparison in Lemma~4.4 and the rational
  scalar cancellation in Lemma~4.5.
  \item The logical reconstruction of Proposition~4.6, Corollary~4.7,
  Proposition~4.9, Theorem~4.10, Theorem~2.3 as used in Section~4, and
  Corollary~4.11 from explicitly separated semantic interfaces.
  \item Type-$D$ norm-gap inequalities and explicit Weyl-orbit packaging for
  Lemma~5.9.
  \item The zero-orbit energy arithmetic used in Proposition~5.4.
  \item Signed-permutation rigidity together with the full type-$D$ root-lattice
  parity argument used in Proposition~5.5.
  \item Appendix-A node-weight derivations used in Proposition~5.11.
  \item Node numerical rigidity used in Proposition~5.12.
  \item The logical combination of the minimal-orbit regularity and rigidity
  steps leading to Propositions~5.10 and~5.13.
  \item The branchwise exhaustion argument leading to Theorem~5.14.
  \item The deductions giving Theorems~6.1 and~6.2.
  \item The additive-subgroup/common-image argument used to produce an actual
  additive-group isomorphism in Theorem~6.3.
  \item The affine-mark arithmetic and coefficient substitution used to pass
  from the abstract BKK/reindexing inputs to Corollary~6.4.
  \item The final assembly of the two conclusions of Theorem~1.1.
\end{enumerate}

\subsection{Representation-theoretic boundary}
The formalization deliberately leaves as explicit inputs results and
infrastructure invoked by the mathematical paper but not reconstructed in
this file.  The boundary includes, among other items, inputs from
Xu/Duflo--Joseph theory, filtered Zhu and Poisson geometry, Li's
$\Delta$-operator and Ramond Zhu theory, finite-$W$/BRST/Skryabin/Premet
machinery, Mili\v{c}i\'c--Soergel and Chen annihilator results, DLM Zhu
correspondence, Petukhov's primitive-ideal/finite-$W$ correspondence,
regular-integral Kazhdan--Lusztig and subregular-cell theory, affine linkage,
and the external results denoted in the source by Jin26, SYZ, and BKK.

These inputs are encoded using structured data such as
\code{PaperData}, \code{MembershipData}, \code{FilteredZhuData},
\code{Section5ProofData}, \code{SYZData}, and
\code{CharacterFormulaData}.  No fake instance is introduced asserting that
\mathlib already contains the actual vertex algebra $L_{-1}(D_\ell)$, its
category $\mathcal O$, the required finite-$W$ modules, BRST functors,
primitive ideals, or the completed affine-Hecke character space.

\section{Detailed audit of the Section 4 membership package}
The structure \code{MembershipData} is the interface used to reconstruct the
membership theorem and simultaneous descent.  A central feature of the current
source is that neither Proposition~4.9 nor Theorem~4.10 is stored as a field.
Instead, the Premet--Whittaker--Skryabin--Chen chain is exposed at lower levels,
and Lean assembles the annihilator identification and ideal containment from
those inputs.

The same four-level convention used below for Section~5 applies here: Level~A
means a concrete calculation or logical deduction proved internally; Level~B
means a paper formula supplied as input with its downstream calculation checked
in Lean; Level~C is a literature-level representation-theoretic theorem exposed
semantically; and Level~D is an encoding bridge between the abstract paper
objects and the concrete formal objects used by the calculation.

\subsection{What is not reconstructed from foundations}
The early twisting construction is not formalized from the definitions of a
minimal $W$-algebra or Li's $\Delta$-operator.  In particular, Lemmas~4.1 and
4.2 are not rederived as standalone theorems about an actual twisted vertex
algebra object.  Their needed consequences are represented through the semantic
setup of the vacuum, cyclic bottom module, charge data, and creation properties.
Thus the existence of the relevant Ramond module and its lowest-space
interpretation remains on the representation-theoretic boundary.

This is a deliberate boundary rather than a hidden axiom: the types
\code{FiniteWModule}, \code{WhittakerModule}, and the predicates describing
quotient modules, highest vectors, simple quotients, and annihilators are all
parameters of \code{MembershipData}.  The file does not claim that \mathlib
already contains the corresponding vertex-algebraic or finite-$W$ objects.

\subsection{Lemmas 4.4--4.5 and Proposition 4.6}
The highest-pair calculation is split into a formula boundary and a concrete
calculation.  The field \code{vacuum\_weight\_from\_li\_twist} supplies equation
(4.5), while Lean proves the Appendix-A coordinate identity
\codeurl{lemma_4_4_premet_shift_from_appendix}, identifying the root-sum shift
with $(r-1)\varpi+\eta_1$.  Hence
\codeurl{lemma_4_4_vacuum_weight} is a Level-B/D-to-A comparison rather than a
direct assumption of Lemma~4.4.

Likewise, \codeurl{central_scalar_from_ramond_relation} supplies the
representation-theoretic scalar formula coming from (4.3) and (4.7), but the
cancellation
\[
-r(r+2)+2(2r+1)\frac{r+1}{4}+\frac{r-1}{2}=0
\]
is proved by \codeurl{lemma_4_5_central_vacuum_zero}.  The zero value in
Lemma~4.5 is therefore Level~A once those formulas are supplied.

For the highest-vector part of Proposition~4.6, the charge ranges
$q\in\{0,1\}$ for current generators and
$q\in\{\tfrac12,\tfrac32\}$ for the chosen odd positive half are exposed as
Level-B/D data, together with the ordinary vacuum creation property.  Lean
then proves in \codeurl{positive_generators_annihilate_vacuum} that all such
positive generators annihilate the vacuum.  The semantic fact that these
families generate Premet's positive part is a Level-C/D bridge.  Combining
this with the internally checked highest-pair calculations gives
\codeurl{proposition_4_6_vacuum_highest_pair}.

\subsection{Corollary 4.7}
The current source no longer assumes directly that the distinguished module
$M^*=L_H(\bar\delta,0)$ is a module over the candidate finite-dimensional
quotient.  Instead, \code{MembershipData} supplies the existence of a simple
quotient of the cyclic bottom module $C^*$, survival of the vacuum highest
line, and Premet's uniqueness of the simple highest-weight quotient.  These
are Level-C/D inputs.  Lean then identifies that simple quotient with
\code{mstar} and transports the quotient-module property, proving
\codeurl{corollary_4_7_mstar_is_quotient_module}.  Thus Corollary~4.7 is a
logical consequence in Lean, although its finite-$W$ highest-weight theory is
external.

\subsection{Lemma 4.8: formalized arithmetic versus semantic comparison}
The rank-one numerical core of Lemma~4.8 is one of the strongest Level-A
pieces in Section~4.  The file defines
\[
c_\beta(m)=\tfrac12 m(m+2),\qquad s_\beta\!\cdot m=-m-2,
\]
and proves \codeurl{cBeta_dot_reflection}, the orthogonal-decomposition
identities corresponding to (4.17)--(4.18), invariance of Premet's second
parameter on the dot orbit, and the specialization at $m=0$.  It also proves
that $s_\beta\cdot0=-2$ is integral and $s_\beta$-antidominant.

The representation-theoretic content of Lemma~4.8 remains Level~C: Premet's
Theorem~1.3, exact Whittaker--Skryabin transport, the comparison of Premet's
parabolically induced module with the Mili\v{c}i\'c--Soergel standard module,
and the uniqueness/transport of simple quotients are fields of
\code{MembershipData}.  Lean combines these pieces at $m=0$ in
\codeurl{lemma_4_8_premet_ms_at_zero}.  Separately,
\codeurl{lemma_4_8_ms_dot_reflection_at_zero} derives the special
$L(0,\psi_\beta)\cong L(-\beta,\psi_\beta)$ statement from the general
Mili\v{c}i\'c--Soergel criterion plus the internally proved equality of
$c_\beta$ parameters.

\subsection{Proposition 4.9 and Theorem 4.10}
Proposition~4.9 is not a trusted conclusion.  The interface provides
annihilator compatibility under the Skryabin/Whittaker model, invariance under
inner conjugation, equality of annihilators for isomorphic Whittaker modules,
Chen's annihilator theorem at an integral antidominant parameter, and the
identification of the highest-weight annihilator at $-\beta$ with $J_{2,\ell}$.
The last Chen hypotheses at $m=-2$ are checked numerically in Lean.
Consequently \codeurl{proposition_4_9_identification} reconstructs
\[
\operatorname{Ann}\operatorname{Skr}(M^*)=J_{2,\ell}
\]
by an explicit chain of equalities.  No field of \code{MembershipData} states
this equality.

The only new literature-level input needed for Theorem~4.10 is the quotient
Skryabin implication saying that a finite-$W$ module over the candidate
quotient has Skryabin annihilator containing $I^{\mathrm{cand}}_\ell$.  This is
encoded as \code{quotient\_skryabin\_containment}.  Lean applies it to the
module obtained in Corollary~4.7 and rewrites the annihilator using the
reconstructed Proposition~4.9.  The result
\codeurl{theorem_4_10_membership} is therefore a genuine derived theorem:
\[
I^{\mathrm{cand}}_\ell\subset J_{2,\ell}.
\]

\subsection{Theorem 2.3 and Corollary 4.11}
The simultaneous-descent step also avoids taking Theorem~2.3 as a single
black box.  The interface separately exposes Xu's reduced-path left-cell
statement, the Duflo--Joseph implication from left-cell equivalence to equality
of primitive annihilators, the model identifying the path weight $\mu_i$ with
the corresponding Weyl element, and the $s_2$ model for $J_{2,\ell}$.  Lean
assembles these in \codeurl{theorem_2_3_path_annihilator}.

Finally, Zhu descent itself remains a Level-C semantic implication.  Combining
that implication with the reconstructed Theorem~2.3 and Theorem~4.10 gives
\codeurl{corollary_4_11_simultaneous_descent}; adding the tautological vacuum
descent produces \codeurl{candidateDescent_from_section4}, exactly the
interface consumed later by Theorem~6.1.

The resulting Section~4 status can therefore be summarized as follows:
\begin{table}[H]
\centering
\scriptsize
\renewcommand{\arraystretch}{1.17}
\begin{tabular}{>{\raggedright\arraybackslash}p{0.20\textwidth}
                >{\raggedright\arraybackslash}p{0.11\textwidth}
                >{\raggedright\arraybackslash}p{0.60\textwidth}}
\toprule
Paper step & Dominant level & Status in the Lean bundle \\
\midrule
Lemmas 4.1--4.2 & C/D & Twisted-module existence/lower-boundedness is not rebuilt from VOA definitions; needed consequences are semantic setup data. \\
Lemma 4.4 & A+B/D & Li-twist weight formula is input; Appendix-A root-coordinate comparison is proved internally. \\
Lemma 4.5 & A+B & Ramond/Premet scalar formula is input; final rational cancellation to zero is proved internally. \\
Proposition 4.6 & A+C/D & Charge inequalities and pair arithmetic are checked; positive-part generation and semantic highest-vector interpretation are boundary inputs. \\
Corollary 4.7 & C/D + internal logic & Simple quotient existence and Premet uniqueness are inputs; identification with $M^*$ and quotient-module conclusion are derived. \\
Lemma 4.8 & A+C & Rank-one parameter algebra is internal; Premet/MS/Skryabin comparison theorems are semantic inputs; the $m=0$ specialization is assembled in Lean. \\
Proposition 4.9 & C/D + internal equality chain & Chen and equivalence/annihilator compatibilities are inputs; $\operatorname{Ann}\operatorname{Skr}(M^*)=J_2$ is derived, not assumed. \\
Theorem 4.10 & C + internal logic & Quotient-Skryabin containment is input; $I^{\mathrm{cand}}\subset J_2$ is derived using Proposition~4.9. \\
Theorem 2.3 / Cor. 4.11 & C/D + internal logic & Xu, Duflo--Joseph and Zhu descent are inputs; common path annihilator and simultaneous descent are assembled in Lean. \\
\bottomrule
\end{tabular}
\caption{Audit of the Section~4 membership/descent chain.}
\end{table}

\section{Detailed audit of the Section 5 proof package}
The structure \code{Section5ProofData} is the main interface through which
Theorem~5.14 is reconstructed.  Importantly, it does not contain
\code{WeightExhaustion} or Theorem~5.14 as a field.  Instead, it exposes the
lower-level data needed for the two associated-variety branches and Lean
assembles the exhaustion theorem from those inputs.

For clarity, the fields are classified into four levels:
\begin{description}[style=nextline,leftmargin=2.4em]
  \item[Level A: internal.] Concrete arithmetic, lattice, Weyl-group, and
  logical deductions are proved in Lean from lower-level data.
  \item[Level B: formula boundary.] A representation-theoretic formula or
  inequality from the paper is supplied, while its substantial downstream
  numerical consequence is proved in Lean.
  \item[Level C: theorem boundary.] A literature-level representation-theoretic
  theorem is supplied through a semantic interface.
  \item[Level D: encoding bridge.] A field connects the abstract predicates in
  \code{PaperData} to concrete coordinates, witnesses, modules, ideals, or
  lattice objects used by the internal calculation.
\end{description}

\subsection{Filtered-Zhu dichotomy}
The field \code{filteredZhu} is a structured Level-C/D interface for
Lemma~3.4.  It supplies the standard Poisson surjection, the PBW--Zhu
filtration identification, the identification of the $C_2$ spectrum with the
minimal-orbit closure, primitive-ideal geometry, Joseph irreducibility, and
the minimal-nonzero-orbit statement.  Lean does \emph{not} take the final
zero/minimal dichotomy as a field.  It proves
\codeurl{poisson_spectrum_containment_from_filtered_zhu}, then
\codeurl{primitive_variety_zero_or_minimal}, and finally
\codeurl{associatedVarietyDichotomy_from_filteredZhu}.  Thus the geometric
inputs are external, but the exact branch dichotomy consumed by
Theorem~5.14 is internally assembled.

\subsection{Zero-orbit branch}
The zero-orbit data consist of \code{BZ}, \code{hcoords}, \code{henergy},
\code{hwitness}, and \code{hroot}.  Their roles are as follows.
\begin{table}[H]
\centering
\scriptsize
\renewcommand{\arraystretch}{1.18}
\begin{tabular}{>{\raggedright\arraybackslash}p{0.20\textwidth}
                >{\raggedright\arraybackslash}p{0.10\textwidth}
                >{\raggedright\arraybackslash}p{0.61\textwidth}}
\toprule
Field & Level & Meaning and checked consequence \\
\midrule
\code{BZ} & D & Supplies $r=\ell-2$, coordinates $(x_1,x_2,c)$, a signed-permutation witness, a root-lattice displacement, and the final identification $(0,0,0)\mapsto\mu_0$. \\
\code{hcoords} & C/D & Encodes the restricted finite-type consequences used in Lemma~5.1: dominance/order, $c\in\{0,1\}$, and $a=x_1-x_2\le r-1$. \\
\code{henergy} & B & Supplies inequality (5.10).  Lean proves from it the strict energy gap $x_1+x_2<2r+1$; the final estimate of Proposition~5.4 is therefore not stored directly. \\
\code{hwitness} & D & Converts the abstract affine-congruence predicate into the explicit signed-permutation/lattice witness used by the concrete type-$D$ calculation. \\
\code{hroot} & D & Supplies the even-sum $D$-root-lattice parity condition for that explicit witness. \\
\bottomrule
\end{tabular}
\caption{Audit of the zero-orbit fields in \code{Section5ProofData}.}
\end{table}

From these fields, the substantive signed-permutation elimination of
Proposition~5.5 is Level~A: the residue analysis and the final parity
contradiction are checked by Lean.  Hence the zero-orbit conclusion
$\mu=\mu_0$ is derived rather than assumed.

\subsection{Minimal-orbit branch}
The field \code{BM} is a \code{MinimalOrbitBridge}.  It separates the
representation-theoretic inputs to Propositions~5.8--5.10 from the numerical
and lattice deductions.  In particular:
\begin{itemize}
  \item \codeurl{q_in_root_lattice} is a Level-D bridge from affine congruence to
  the explicit $D$-root-lattice displacement;
  \item \code{centralGapData} contains Level-B/C constituent data together
  with the original weighted-trace identity (3.5) and equation (4.7);
  \item \code{lemma57IntegralityData} is a Level-C/D package for the chain
  ``affine integrality $\Rightarrow$ Lemma~5.6 highest pair $\Rightarrow$
  cyclic root/vector shifts'' used in Lemma~5.7;
  \item \codeurl{normGap_eq_central_scalar} is the Level-B identity (5.13); and
  \item \codeurl{petukhov_node_classification} is a Level-C instance of
  Petukhov's classification of the regular minimal primitive fibre.
\end{itemize}

On top of this boundary, Lean derives several nontrivial steps internally.
The four-way level-one $D_r$ weight classification, exclusion of spinor
constituents from integral coordinates, the constituentwise Casimir bound,
and the dimension-weighted estimate
$C(M_J)\le C_{\max}$ are proved in Lean.  The algebraic elimination of the
Ramond conformal scalar from (3.5) and (4.7) is also internal, yielding
$C=2C(M_J)-r(r+2)$.  Combined with the fully formalized type-$D$ norm gap of
Lemma~5.9, Lean proves Proposition~5.10 in the regular-central-character
form used later.

\subsection{Node-2 rigidity and the Proposition 5.11 subtlety}
The field \code{hnodeBounds} returns a \code{PrimitiveNodeBoundsBridge}.
It contains three semantic links: candidate-ideal containment makes the
node module a $B_\ell$-module; Lemma~3.1 supplies the $A_1$ and $D_r$ bounds;
and the actual node module's invariants agree with the numerical values of
Proposition~5.11.

The last item must be distinguished from the coordinate calculation itself.
The appendix formulas and the resulting node weights are independently
proved in Lean, including the $D_5$ boundary.  What remains at Level~D is the
identification of those formal coordinate expressions with the invariants of
the semantically represented finite-$W$ node modules.  Once this bridge is
supplied, the elimination of every node except node~2 is Level~A and is
proved by \codeurl{nodeAdmissible_forces_value_one} and the integrated
Proposition~5.12 theorem.

\subsection{Final left-cell step}
The field \code{leftCellBridge} exposes rather than hides the two external
ingredients in the end of Proposition~5.13: regular-integral
primitive-ideal equality implies the Kazhdan--Lusztig left-cell relation,
and the finite subregular left cell of $s_2$ consists of the path words.
These are Level-C inputs.  The final conversion of a path word into one of
the weights $\mu_i$ is represented separately, and Lean combines regularity,
Petukhov's node classification, node-2 rigidity, and the left-cell bridge to
obtain the minimal-orbit path-weight conclusion.

Consequently, \codeurl{weightExhaustion_from_integrated_theorem_5_14} is a
genuine Lean theorem: the final exhaustion statement is not present in
\code{Section5ProofData}.  What is trusted is the explicitly listed
representation-theoretic boundary beneath it.

\section{Detailed audit of Section 6 and the top-level theorem}
Section~6 is where the outputs of the membership and exhaustion arguments are
converted into the stated simple-object classification and then into the
Grothendieck-group comparison.  The current Lean bundle exposes this passage
at a relatively fine level.  In particular, Theorem~6.1, Theorem~6.2, and the
additive-group content of Theorem~6.3 are not stored as direct result fields of
the top-level data structure.  They are assembled from lower-level interfaces.

As above, Level~A denotes an internal Lean proof, Level~B a paper formula or
structural identity supplied as input with its consequences checked in Lean,
Level~C a literature-level representation-theoretic theorem, and Level~D an
encoding bridge connecting the semantic paper objects to the abstract types
used by the formalization.

\subsection{Theorem 6.1: classification for the candidate quotient}
The standalone structure
\codeurl{DTypeCandidateQuotientClassification.ClassificationData} contains
only the underlying weight and module types, the vacuum and path labels, the
module constructor, and the two predicates \code{InAmbientBlock} and
\code{FactorsThroughQ}.  The theorem
\codeurl{DTypeCandidateQuotientClassification.ClassificationData.theorem_6_1}
then proves
\[
  (\text{ambient block and factors through }Q_\ell)
  \quad\Longleftrightarrow\quad
  \text{vacuum or one of the path modules}.
\]
Its reverse implication is a purely logical use of the candidate ambient-block
and descent data.  Its forward implication first uses the affine
highest-weight parametrization and then invokes the exhaustion statement for
the extracted finite highest weight.

The top-level theorem does not assume the Section~4 descent conclusion or the
Section~5 exhaustion conclusion directly.  The field \code{S4} is converted
by \codeurl{candidateDescent_from_section4} into \code{CandidateDescent}, and
\code{S5} is converted by
\codeurl{weightExhaustion_from_integrated_theorem_5_14} into
\code{WeightExhaustion}.  Lean then combines this with
\codeurl{FactoringWeightConditions} to construct the exact
\code{BlockExhaustion} interface needed by Theorem~6.1.

The remaining Section~6 boundary at this stage is explicit.  The hypothesis
\codeurl{CandidateInAmbientBlock} records that the vacuum and path modules belong
to the ambient block (Level~C/D).  The hypothesis \codeurl{AmbientParametrization}
encodes affine linkage/highest-weight parametrization of every ambient-block
simple (Level~C).  Finally, \codeurl{FactoringWeightConditions} packages the
representation-theoretic passage from a factoring module to the two finite
weight conditions to which Theorem~5.14 applies: the candidate Zhu-ideal
condition and the affine congruence (Level~C/D).  Once these are supplied,
Theorem~6.1 itself is Level~A.

\subsection{Theorem 6.2 and the statement ``exactly $\ell+1$ simples''}
The structure \codeurl{DTypeSimpleObjects.SimpleBlockData} separates the
candidate quotient predicate from the predicate \code{IsLminusOneModule}.
The sole new representation-theoretic bridge needed for Theorem~6.2 is
\codeurl{QuotientIsomorphismInput}, which records, at the classification level,
that an ambient-block simple factors through $Q_\ell$ if and only if it is a
module for $L_{-1}(D_\ell)$.  This is the formal interface for the external
isomorphism
\[
  Q_\ell\cong L_{-1}(D_\ell)
\]
and is Level~C.  Given Theorem~6.1 and this bridge,
\codeurl{theorem_6_2_from_theorem_6_1} proves the classification equivalence
internally.

The numerical phrase ``exactly $\ell+1$ simple objects'' is formalized more
carefully than a bare cardinality assertion.  Lean uses the label type
\[
  \code{Option (Fin ell)},
\]
where \code{none} is the vacuum and \code{some i} is a path module.  The theorem
\codeurl{label_cardinality} proves internally that this label type has
cardinality $\ell+1$.  The theorem \codeurl{simple_objects_equiv_labels} then
constructs an actual equivalence between this label type and the subtype of
simple vacuum-block $L_{-1}(D_\ell)$-modules.

Two semantic points remain visible in that construction.  First,
\codeurl{PaperData.moduleOfWeight_is_simple} declares that
\code{moduleOfWeight mu} denotes the simple affine highest-weight module with
that highest weight; simplicity is not rebuilt from an affine Kac--Moody
module definition.  Second, the top-level hypothesis
\codeurl{CandidateLabelsDistinct} asserts that the vacuum and path modules are
pairwise distinct.  The current bundle does not derive this distinctness from
an internal comparison of the explicit highest weights.  Thus the bijection
construction and the cardinality of the abstract label set are Level~A, while
the semantic simplicity and distinctness inputs are Level~D/C.

\subsection{Theorem 6.3: common ambient image and additive equivalence}
The structure \codeurl{DTypePaperSYZFinal.SYZData} is the semantic boundary
for the Shan--Yan--Zhao comparison.  It provides abstract additive groups
\code{K0}, \code{Cell}, and \code{Ambient}; injective additive maps
\code{embedK0} and \code{embedCell}; indexed simple classes and cell
basis generators; and their designated ambient basis elements.  The existence
and injectivity of these canonical ambient realizations are Level~C/D inputs.

The four hypotheses passed to \codeurl{theorem_6_3_isomorphism} are deliberately
primitive:
\begin{itemize}
  \item the simple classes generate the encoded Grothendieck group;
  \item the cell generators generate the specialized cell module;
  \item the $K_0$ generators map to the designated ambient basis elements;
  \item the cell generators map to the same ambient basis elements.
\end{itemize}
These are representation-theoretic/categorical boundary facts (Level~C/D),
not conclusions silently built into \code{SYZData}.

From them, the generic theorem
\codeurl{range_eq_CanonicalSpan_of_generated} proves internally that each
canonical embedding has range equal to the same additive subgroup generated by
the indexed ambient basis.  Consequently
\codeurl{theorem_6_3_final} proves equality of the two ranges.  Lean then
constructs explicit additive equivalences from $K_0$ and from the cell module
to this common image and composes them to obtain
\codeurl{commonImageFactorizationAddEquiv}.  Hence
\codeurl{theorem_6_3_isomorphism} produces an actual additive-group
isomorphism rather than merely an equality of cardinalities or a set-level
bijection.  The separate theorem
\codeurl{theorem_6_3_ambient_compatibility} verifies that this equivalence
intertwines the two ambient embeddings.

There is an important scope limitation.  The formal theorem is an isomorphism
of additive groups.  The source does not construct the affine-Hecke or affine
Weyl action on these groups and does not prove an equivariance statement for
that action.  Thus the kernel-checked content matches the common-image
additive-group comparison, not a foundational formalization of the full module
structure used in the surrounding representation theory.

\subsection{Corollary 6.4: uniform character formula}
The uniform character formula is formalized separately from the top-level
Theorem~1.1 wrapper.  The concrete affine type-$D$ mark function
\codeurl{affineMark} and its endpoint/interior properties are proved in Lean.
The structure \codeurl{CharacterFormulaData}, however, deliberately leaves the
completed character space, the finite and lattice sums, the BKK affine-Weyl
expansion, the paired reindexing, and the BKK multiplicity-difference identity
as external data.  These are Level~B/C interfaces, because constructing the
completed formal character space and the BKK theory is outside the bundle.

Once those three character-theoretic identities are supplied, Lean proves that
the multiplicity difference is exactly the displayed coefficient
$b_i(u,\gamma)$, substitutes that identity under the root-lattice and finite
Weyl sums, and derives
\codeurl{corollary_6_4_uniform_characters}.  Accordingly, the coefficient
substitution and final formal-sum assembly are Level~A, while the BKK expansion
and affine-Weyl pairing/reindexing remain external.  Corollary~6.4 is not one
of the outputs bundled into \codeurl{DTypeMainTheorem.theorem_1_1}; it is a
separate verified consequence in the source.

\subsection{The final wrapper and its exact trust boundary}
The theorem \codeurl{DTypeMainTheorem.theorem_1_1} visibly lists the remaining
semantic inputs in its signature.  The most important are
\codeurl{CandidateInAmbientBlock}, \codeurl{AmbientParametrization},
\codeurl{FactoringWeightConditions}, \codeurl{QuotientIsomorphismInput},
\codeurl{CandidateLabelsDistinct}, and the Shan--Yan--Zhao ambient realization and
generator-map data.  By contrast, the Section~4 descent theorem and the
Section~5 exhaustion theorem are not top-level assumptions: they are derived
inside the proof from \code{S4} and \code{S5}.

The wrapper then performs four kernel-checked assembly steps.  It derives the
candidate-quotient classification, transports it across
$Q_\ell\cong L_{-1}(D_\ell)$, constructs the label equivalence and proves the
$\ell+1$ label cardinality, and constructs the additive $K_0$--cell-module
isomorphism.  This gives a precise reading of the final verification claim:
Lean checks the global deductive architecture from the explicitly listed
representation-theoretic interfaces, while the mathematical validity of those
interfaces remains delegated to the paper and the cited literature.

\begin{table}[H]
\centering
\scriptsize
\renewcommand{\arraystretch}{1.18}
\begin{tabular}{>{\raggedright\arraybackslash}p{0.17\textwidth}
                >{\raggedright\arraybackslash}p{0.12\textwidth}
                >{\raggedright\arraybackslash}p{0.62\textwidth}}
\toprule
Section-6 item & Dominant level & Status in the Lean bundle \\
\midrule
Theorem 6.1 & A+C/D & Descent and exhaustion are reconstructed from Sections~4--5; ambient parametrization and factoring-weight bridges remain external. \\
Theorem 6.2 & A+C & Classification is transported internally across the external $Q_\ell\cong L_{-1}(D_\ell)$ interface. \\
Exactly $\ell+1$ simples & A+D/C & Label cardinality and bijection construction are internal; semantic simplicity and pairwise distinctness of candidate modules are boundary inputs. \\
Theorem 6.3 & A+C/D & Common-span/range equality and an actual additive equivalence are constructed internally from injective ambient realizations, generation, and generator-image data. \\
Corollary 6.4 & A+B/C & Type-$D$ marks and coefficient substitution are internal; BKK expansion, paired reindexing, and completion are external. \\
Theorem 1.1 wrapper & A+C/D & Final classification, count, and additive isomorphism are assembled in Lean from the explicit interfaces in the theorem signature. \\
\bottomrule
\end{tabular}
\caption{Boundary audit for Section~6 and the final wrapper.}
\label{tab:section6audit}
\end{table}

\section{Paper-to-Lean verification map}
Table~\ref{tab:map} records the intended relationship between representative
results of the paper and the Lean bundle.

\begin{table}[H]
\centering
\scriptsize
\renewcommand{\arraystretch}{1.22}
\begin{tabular}{>{\raggedright\arraybackslash}p{0.17\textwidth}
                >{\raggedright\arraybackslash}p{0.32\textwidth}
                >{\raggedright\arraybackslash}p{0.43\textwidth}}
\toprule
Paper result & Lean declaration / status & Verification boundary \\
\midrule
Lemma 4.4 & \codeurl{lemma_4_4_vacuum_weight} & Li-twist weight formula is a bridge input; the Appendix-A root-coordinate identity is kernel-checked. \\
Lemma 4.5 & \codeurl{lemma_4_5_vacuum_central_zero} & Final scalar cancellation is kernel-checked from the Ramond/Premet formula boundary. \\
Proposition 4.6 & \codeurl{proposition_4_6_vacuum_highest_pair} & Positive-generator annihilation and pair assembly are derived from explicit charge/generator interfaces. \\
Corollary 4.7 & \codeurl{corollary_4_7_mstar_is_quotient_module} & Derived from simple-quotient and Premet-uniqueness interfaces. \\
Lemma 4.8 & rank-one concrete lemmas and \codeurl{lemma_4_8_premet_ms_at_zero} & Parameter algebra is internal; Premet/MS/Skryabin comparison remains semantic. \\
Proposition 4.9 & \codeurl{proposition_4_9_identification} & Annihilator equality is reconstructed; it is not a field of \code{MembershipData}. \\
Theorem 4.10 & \codeurl{theorem_4_10_membership} & Ideal containment is derived from quotient-Skryabin input and reconstructed Proposition~4.9. \\
Corollary 4.11 & \codeurl{corollary_4_11_simultaneous_descent} & Theorem~2.3 and Zhu descent are assembled with Theorem~4.10. \\
Proposition 5.4 & \codeurl{proposition_5_4_energy_gap} & Arithmetic inequality kernel-checked. \\
Proposition 5.5 & \codeurl{proposition_5_5_full_root_lattice} and integration lemmas & Signed permutations and type-$D$ parity encoded explicitly. \\
Lemma 5.9 & norm-gap/Weyl-orbit modules and integration bridge & Integer inequalities, root-lattice parity and orbit packaging checked internally. \\
Proposition 5.10 & \codeurl{proposition_5_10_from_5_8_and_5_9} & Deduction is kernel-checked from Proposition 5.8 boundary data and the internal norm gap. \\
Proposition 5.11 & \codeurl{proposition_5_11_from_appendix} & Type-$D$ node-weight calculation checked internally. \\
Proposition 5.12 & \codeurl{proposition_5_12_integrated} & Numerical rigidity checked internally from explicit hypotheses. \\
Proposition 5.13 & \codeurl{proposition_5_13} / integrated path conclusion & Deduction checked from explicit primitive/KL and node data. \\
Theorem 5.14 & \codeurl{theorem_5_14_from_propositions_5_5_and_5_13} & Zero- and minimal-orbit branches assembled in Lean. \\
Theorem 6.1 & \codeurl{theorem_6_1} & Descent/exhaustion are internally supplied from Sections~4--5; ambient parametrization and factoring-weight bridges remain external. \\
Theorem 6.2 & \codeurl{theorem_6_2_from_theorem_6_1} & Classification is deduced in Lean from Theorem~6.1 and the external $Q_\ell\cong L_{-1}(D_\ell)$ interface. \\
Theorem 6.3 & \codeurl{theorem_6_3_isomorphism} & Common-range equality and an actual additive-group isomorphism are constructed from generator and ambient-image data; no affine-action equivariance is claimed. \\
Corollary 6.4 & \codeurl{corollary_6_4_uniform_characters} & Coefficient substitution and final formal-sum assembly are internal from external BKK/reindexing inputs. \\
Theorem 1.1 & \codeurl{DTypeMainTheorem.theorem_1_1} & Final classification, $\ell+1$ label count, and additive-group isomorphism assembled from all explicit inputs. \\
\bottomrule
\end{tabular}
\caption{Representative paper-to-Lean verification map.}
\label{tab:map}
\end{table}

\section{The top-level theorem}
The theorem \codeurl{DTypeMainTheorem.theorem_1_1} assumes $5\le \ell$ and a
structured package \code{PaperData ell}.  It receives explicit Section~4
membership data, Section~5 representation-theoretic proof data, ambient
parametrization and factoring-weight bridges, the quotient-isomorphism input,
distinctness of candidate labels, and the ambient Shan--Yan--Zhao realization.
Notably, neither \code{CandidateDescent} nor \code{WeightExhaustion} is a
top-level hypothesis: both are derived inside the proof from the Section~4 and
Section~5 packages.  From these inputs, Lean derives three outputs:
\begin{enumerate}
  \item the classification predicate corresponding to equation (1.1);
  \item an equivalence between the finite candidate-label type and the simple
  vacuum-block objects, together with a proof that the label type has
  cardinality $\ell+1$; and
  \item a nonempty additive-group isomorphism between the encoded
  Grothendieck group and the encoded specialized cell module, corresponding
  to equation (1.2).
\end{enumerate}
This formulation makes the trust boundary visible in the theorem signature.

\section{Source-level audit}
The archived source in this version incorporates the latest proof-term closure
fix in \codeurl{middleLevel_from_coefficients}: after rewriting the middle-mass
identity and pushing the cast, the final equality is closed explicitly by
\code{rfl}.  This change does not alter the representation-theoretic boundary
or the top-level theorem signature; it makes the local coefficient calculation
syntactically closed in the released source.

For the archived source file used for this report, the following static facts
were checked directly from the file:
\begin{itemize}
  \item the source has 10,511 lines;
  \item it contains 304 actual declarations introduced with \code{theorem} or
  \code{lemma} (261 theorems and 43 lemmas);
  \item the source contains no proof placeholder introduced by \code{sorry}
  or \code{admit};
  \item it contains no user-declared Lean \code{axiom}; and
  \item the final source comment records the formalization boundary and
  distinguishes the kernel-checked deduction from a foundational construction
  of the underlying representation-theoretic objects.
\end{itemize}
The SHA-256 digest of the archived Lean source is
\begin{center}
{\scriptsize\ttfamily 8eeff0283a801501d18c4f70c83b970d4ab8139ca7acd7c9262ae04b1848f701}.
\end{center}

This source-level statement should be distinguished from a dependency audit of
Lean's logical constants.  Under the frozen environment, the command
\codeurl{\#print axioms DTypeMainTheorem.theorem_1_1} reports exactly
\begin{center}
\code{[propext, Classical.choice, Quot.sound]}.
\end{center}
Thus the final theorem uses only these standard foundational Lean axioms; no
user-declared representation-theoretic axiom occurs in the final proof term.

\section{Reproducibility}
The ancillary release contains the Lean source together with the project
metadata supplied for the verification environment.  The frozen toolchain is
\code{leanprover/lean4:v4.34.0-rc1}.  The Lake project requires Mathlib at tag
\code{v4.34.0-rc1}, while \code{lake-manifest.json} resolves Mathlib to the
exact commit
\begin{center}
{\scriptsize\ttfamily de5ce8a9a66a4aa68a9bdbb35b63a06d34d9ca11}.
\end{center}
The ancillary directory therefore includes
\codeurl{lean-toolchain}, \codeurl{lakefile.toml}, and
\codeurl{lake-manifest.json} in addition to the Lean source.

From the project directory, the source is checked directly by
\begin{center}
\code{lake env lean SubregularAffineCells.lean}.
\end{center}
A verification run on the released working source completed with process exit
code~0.  The executable identified itself as
\code{Lean 4.34.0-rc1} for \codeurl{x86_64-w64-windows-gnu}.  Its Lean commit was
\begin{center}
{\scriptsize\ttfamily 3447a668783dbce1a8fdb97101dd067687b2b418}.
\end{center}
With the temporary command
\begin{center}
\codeurl{\#print axioms DTypeMainTheorem.theorem_1_1}
\end{center}
appended after the final namespace, Lean reported
\begin{center}
\code{[propext, Classical.choice, Quot.sound]}.
\end{center}
The frozen manifest fixes the dependency revisions, while this axiom audit
records the foundational constants occurring in the final proof term.  These
are separate from the absence of user-declared representation-theoretic axioms
in the source.

\section{What this verification does and does not claim}
The strongest accurate summary of the formalization is:
\begin{quote}
The Lean bundle provides a kernel-checked internal deduction of the principal
paper conclusions from explicit representation-theoretic boundary inputs.
\end{quote}
It should \emph{not} be described as a foundational formalization from the
definitions of vertex operator algebras, BRST reduction, finite $W$-algebras,
primitive ideals, or affine Hecke theory.  This distinction is essential both
for mathematical accuracy and for reproducibility.

\section{Availability of formalization}
The Lean source accompanying this report is distributed with the arXiv source
submission as the ancillary file
\codeurl{anc/SubregularAffineCells.lean}.  The frozen reproducibility metadata
are provided alongside it as \codeurl{anc/lean-toolchain},
\codeurl{anc/lakefile.toml}, and \codeurl{anc/lake-manifest.json}; the three
boundary-audit notes and \codeurl{anc/VERIFICATION.txt} are included in the
same directory.

A public, maintainable copy is hosted at
\begin{center}
\url{https://github.com/jshemail12345-debug/SubregularAffineCells-Lean}.
\end{center}
The exact GitHub snapshot intended to correspond to the initial arXiv version
of this verification report is the tagged release
\begin{center}
\url{https://github.com/jshemail12345-debug/SubregularAffineCells-Lean/releases/tag/v1.0.0}.
\end{center}
The arXiv ancillary copy and the tagged GitHub release therefore serve
complementary purposes: the former is attached to the report version, while
the latter provides a fixed, clonable public repository snapshot.

\section*{Acknowledgements}
The author acknowledges the assistance of the language model
\textsc{DeepSeek-V4-Pro} in drafting and refining portions of the Lean~4
formalization and in preparing this verification report.  The author remains
responsible for the mathematical statements, the formalization boundary, the
verification claims, and the final contents of the released source and report.

\end{document}